\documentclass [a4paper,twoside,11pt]{article}
\usepackage{amsmath,amsfonts,amssymb}
\usepackage{vmargin,graphicx,theorem}
\usepackage[english]{babel}
\usepackage{enumerate}
\usepackage{color}

\setpapersize[portrait]{A4}
\setmarginsrb{2.5cm}{2.5cm}{2.5cm}{2cm}{0cm}{0.1cm}{0.5cm}{2cm}
\newcommand{\disp}{\displaystyle}
\newcommand{\dR}{\ensuremath{\mathbb{R}}}

\newtheorem{ethm}{Theorem}[section]

\newtheorem{ecor}[ethm]{Corollary}

\newtheorem{erem}[ethm]{Remark}

\newcommand{\proofend}{~$\rhd$}
\newcommand{\proofbegin}{~$\lhd$}

\newenvironment{eproof}
               {\noindent {\emph{\textbf{Proof}}}\\\proofbegin~}
               {\proofend\\}

\newcommand{\p}[4]{{#3}\!\left#1{#4}\right#2}

\newcommand{\ABS}[1]{\ensuremath{{\left| #1 \right|}}} % |1|
\newcommand{\PAR}[1]{\ensuremath{{\left(#1\right)}}} % (1)
\newcommand{\SBRA}[1]{\ensuremath{{\left[#1\right]}}} % [1]
\newcommand{\BRA}[1]{\ensuremath{{\left\{#1\right\}}}} % {1}
\newcommand{\NRM}[1]{\ensuremath{{\left\Vert #1\right\Vert}}} % ||1||
\renewcommand{\phi}{\varphi}
\renewcommand{\geq}{\geqslant}

\newcommand{\varf}[1]{\rm{Var}_{#1}}
\newcommand{\entf}[1]{{\rm{Ent}}_{#1}}

\newcommand{\ent}[2]{\p(){\entf{#1}}{#2}}

\newcommand{\var}[2]{\p(){\varf{#1}}{#2}}

\def\disp{\displaystyle}

\newcommand{\A}{\ensuremath{\mathcal A}}

\newcommand{\GI}{\mathbf{L}}
\newcommand{\PT}[1]{\mathbf{P}_{\!\bf#1}}
\newcommand{\Pt}[1][t]{\ensuremath{\mathbf{P_{\!#1}}}}
\newcommand{\NT}[1]{\mathbf{N}_{\!#1}}
\newcommand{\Nt}[1][t]{\ensuremath{\mathbf{N_{\!\bf#1}}}}
\newcommand{\TT}[1]{\mathbf{T}_{\!\bf#1}}
\newcommand{\Tt}[1][t]{\ensuremath{\mathbf{T_{\!\bf#1}}}}
\newcommand{\LT}[1]{\mathbf{L}_{\bf#1}}
\newcommand{\Lt}[1][t]{\ensuremath{\mathbf{L_{#1}}}}

\newcommand{\QT}[1]{\mathbf{Q}_{\bf#1}}
\newcommand{\Qt}[1][t]{\ensuremath{\mathbf{Q_{#1}}}}

\newcommand{\beq}{\begin{equation}}\newcommand{\eeq}{\end{equation}}
\begin{document}

\title{\sl Dimension dependent  hypercontractivity for Gaussian kernels }
\author{Dominique Bakry\thanks{Universit\'e Paul-Sabatier, Institut de Math\'ematiques de Toulouse}~\thanks{Institut Universitaire de France.} 
 ,  Fran{\c c}Åois Bolley\thanks{Universit\'e Paris-Dauphine, Ceremade.} 
 \,and Ivan Gentil\footnotemark[3] }

\date{\today}

\maketitle

\abstract{We derive sharp, local and dimension  dependent hypercontractive bounds on the Markov kernel of a large class of  diffusion semigroups. Unlike  the dimension free ones, they capture refined properties of Markov kernels, such as  trace estimates. They  imply classical bounds on the Ornstein-Uhlenbeck semigroup and  a dimensional and refined (transportation) Talagrand inequality when applied to the Hamilton-Jacobi equation.  Hypercontractive bounds on the Ornstein-Uhlenbeck semigroup driven by a non-diffusive L\'evy semigroup are also  investigated. Curvature-dimension criteria are the main tool in the analysis.  }

\bigskip


\bigskip

%%%%%%%%%%%%%%%%%%%%%%%%%%%%%%%%%%%%%%%%%%%%%%%%%%%%%%
%%%%%%%%%%%%%%%%%%%%%%%%%%%%%%%%%%%%%%%%%%%%%%%%%%%%%%   
 \section{Introduction}

The well known Nelson Theorem asserts that the Ornstein-Uhlenbeck semigroup $(\Nt)_{t\geq0}$  satisfies the hypercontractive bound  
\begin{equation}
\label{eq-hyper1}
\NRM{\Nt f}_{L^{q_2}(d\gamma)}\leq \NRM{f}_{L^{q_1}(d\gamma)}
\end{equation}
 for all $q_2\geq q_1>1$ such that $q_2-1\leq e^{2t}(q_1-1)$; here $\gamma$ is the standard Gaussian measure on $\dR^n$.  It means that for given $q_2\geq q_1>1$, there exists a large enough time $t$ such that $N_t f\in L^{q_2}(d\gamma)$ if $f\in L^{q_1}(d\gamma)$ only. 
   
 For a general diffusion Markov semigroup, the Gross Theorem (see~\cite{gross}) shows that this hypercontractivity property is equivalent to a logarithmic Sobolev inequality on its invariant measure.  Let indeed $(\Pt)_{t\geq0}$ be such a semigroup, with invariant measure $\mu$ and Dirichlet form $\mathcal E_\mu(f,f)$. We say that $\mu$ satisfies a logarithmic Sobolev inequality with constant $C$ if
 $$
 \ent{\mu}{f^2}:=\int f^2\log \frac{f^2}{\int f^2d\mu}d\mu\leq C\mathcal E_\mu(f,f)
 $$
for all smooth functions $f.$ Then the Gross Theorem asserts that this inequality is equivalent to the hypercontractive bound
$$
\NRM{\Pt f}_{L^{q_2}(d\mu)}\leq \NRM{f}_{L^{q_1}(d\gamma)}
$$ 
for all $q_2\geq q_1>1$ such that $q_2-1\leq e^{4t/C}(q_1-1)$. For instance the Gaussian measure is the invariant measure for the Ornstein-Uhlenbeck semigroup $(\Nt)_{t \geq 0}$. 

In many cases (Euclidean spaces, Riemannian manifolds, etc), the Dirichet form is simply $\int |\nabla f|^2 d\mu$ and it is often convenient to rewrite the logarithmic Sobolev inequality as
 $$
 \int f\log f d\mu -\int f d\mu \; \log \int f d\mu \leq \frac{C}{4} \int \frac{|\nabla f|^2}{f} d\mu
 $$ 
for all positive functions $f$,  the left hand side being an entropy and the right hand side a Fisher information.

These properties of semigroups are restricted to semigroups for which the invariant measure $\mu$ is finite. But there are many interesting semigroups, for example the usual heat semigroup in the Euclidean space, for which the invariant measure is infinite (it is the Lebesgue measure in the case of the heat semigroup).  For these semigroups, one still have logarithmic Sobolev inequalities, but for the heat kernel measure $\Pt$ itself instead of the invariant measure $\mu$. For example, the Ornstein-Uhlenbeck semigroup $(\Nt)_{t \geq 0}$ satisfies     the so-called local logarithmic Sobolev inequality 
\beq\label{ineq:LocLSNt}
 \ent{\Nt}{f}:=\Nt( f\log {f})- \Nt( f)\log \Nt( f)\leq \frac{1 - e^{-2 t}}{2} \, {\Nt}\PAR{\frac{|\nabla f|^2}{f}}
\eeq
for all positive $f$ and $t\geq0$.  The name local comes from the fact that for small time $t$ the heat kernel measure $\Nt(x,dy)$  is concentrated around $x$ as a Gaussian measure of small variance (but it is also local in time).   Letting $t$ go to infinity, then $\Pt f$ goes to $\int fd\gamma$ and one recovers the logarithmic Sobolev inequality  for the ergodic measure $\gamma$.

On the other hand, the Euclidean heat kernel  $\Pt$ on $\dR^n$ also satisfies a local  logarithmic Sobolev inequality 
\beq\label{ineq:LocLSEucl}
\ent{\Pt}{f} :=\Pt (f\log f)-\Pt f\log\Pt f\leq t \, \Pt(\frac{|\nabla f|^2}{f}),
\eeq
 which in this context is less useful since here one cannot let $t$ go to infinity (here the invariant measure, up to a scaling constant, is $\lim_{t\to \infty} t^{n/2}\Pt$).

These local inequalities are in fact very general, and only depend on a curvature bound on the generator, more precisely, a $CD(\rho, \infty)$ inequality, and they are equivalent to it. We shall give more details  later on. They do not depend on the dimension, which can be a great advantage with the idea of easy extensions to infinite dimensional settings when this  makes sense; but it is also a drawback since they do not capture the precise dimensional information about the semigroup. 

For the Euclidean heat kernel $\Pt$, the local logarithmic Sobolev inequality (\ref{ineq:LocLSEucl}) can be improved into the dimensional bound \beq\label{ineq:LSdimEucl}
\ent{\Pt}{f} \leq t \, \Delta \Pt f+\frac{n}{2}\Pt f\log\PAR{1-\frac{2t}{n}\frac{\Pt\PAR{f\Delta (\log f)}}{\Pt f}}
\eeq
for all positive $f$ and $t$. As proved in~\cite{bakry-ledoux06}, this bound is in fact equivalent to the more precise $CD(0,n)$ inequality of the Euclidean space, and is equivalent to it. It improves on~\eqref{ineq:LocLSEucl} since
$$
\Delta(\log f) = \frac{\Delta f}{f}-\frac{|\nabla f|^2}{f^2}
$$
and $\log (1+x)\leq x$ for all $x$.

Inequality~\eqref{ineq:LSdimEucl} also has a reverse form
\beq \label{ineq:LiYauEucl} 
\ent{\Pt}{f} \geq  t\, \Delta \Pt f  -{n\over 2} \,\Pt f
      \log  { } \bigg ( 1 + {2t \over n} \, \Delta   (\log \Pt f)  \bigg ),
  \eeq which is in fact a reinforced form of the celebrated Li-Yau inequality (see \cite{bakry-ledoux06,li-yau}).
  
  These inequalities may still be stated in the more general setting of a diffusion semigroup on a space $E$, for which they appear to be equivalent to a $CD(0, n)$ criterion. For example, for the heat kernel on a Riemannian manifold, this means that they are equivalent to the Ricci curvature being non negative.
  
  One should be careful when extending such equivalences between functional properties of the heat kernel and lower bounds on the Ricci curvature to situations where the generator of the semigroup is not the Laplace Beltrami operator of a Riemannian manifold:  both inequalities (\ref{ineq:LocLSEucl}) and (\ref{ineq:LSdimEucl}) are equivalent to the non negativity of the Ricci curvature for such heat kernels, but it is no longer the case for general diffusion semigroups. In particular, the Ornstein-Uhlenbeck semigroup $(\Nt)_{t \geq 0}$ will never satisfy any dimensional inequality like (\ref{ineq:LSdimEucl}).

\bigskip

In this paper we describe how a local logarithmic Sobolev inequality can be interpreted as a local hypercontractive bound on the semigroup. 
In fact we should be careful that if for instance for fixed $t$ we brutally apply  the Gross Theorem to the local logarithmic Sobolev inequality (\ref{ineq:LocLSEucl}), we obtain a hypercontractivity property not for the semigroup $(\Pt)_{t \geq 0}$ itself, but for the semigroup which has the Markov kernel $\Pt (x,dy)$ as invariant measure and same ``carr\'e du champ". With a few exceptions this semigroup is in general not  easily related to~$\Pt$. 

In sections~\ref{sec-cdrho} and~\ref{sec-cdn} we investigate the consequences on the semigroup itself of local logarithmic Sobolev inequalities such as~\eqref{ineq:LocLSNt} and~\eqref{ineq:LSdimEucl}. We observe that many Markov semigroups satisfy a local and dimensional hypercontractive bound. As an application this gives a new characterization of the curvature-dimension criterion without any derivative.  

In particular, when applied to the Euclidean heat semigroup, they will turn into dimensional hypercontractive bounds on the Ornstein-Uhlenbeck semigroup which make the Nelson theorem more precise (see section~\ref{sec-appliOU}).
As a simple consequence we recover sharp bounds on the Ornstein-Uhlenbeck Markov kernel and on its trace, which brings new links between these quantities and curvature  bounds, and ensure the optimality of our method. 
% deduce from this that $\Nt$ is Hilbert-Schmidt (which of course is not new), by linking this property to curvature bounds (which is new). 

In section~\ref{sec-transport} we draw some consequences in terms of dimensional transportation inequalities, through the Hamilton-Jacobi equations; in particular we obtain a local transportation inequalities for a large class of semigroups and give dimensional improvements on a celebrated transportation inequality for the Gaussian measure by M. Talagrand. 

Non-diffusion situations of Levy semigroups are finally studied in section~\ref{sec-levy}. We prove hypercontractive bounds on the  Ornstein-Uhlenbeck semigroup driven by an $\alpha$-stable L\'evy process :  to our knowledge this is the first example of a hypercontractivity property for a non-diffusive semigroup with non-reversible invariant measure.   

\medskip

The following section is devoted to a presentation of the setting and tools that will be used throughout this work.

%%%%%%%%%%%%%%%%%%%%%%%%%%%%%%%
\section{Diffusion semigroups and $\Gamma_2$ criterion}
\label{sec-diff}
%%%%%%%%%%%%%%%%%%%%%%%%%%%%%%%%

The $\Gamma_2$  calculus and curvature-dimension inequalities are an easy way to get hypercontractivity results, both for the measure $\Pt (x,dy)$ and for the invariant measure $\mu$.  It relies on the local analysis of the generator $L$ of the Markov semigroup. 

 A Markov semigroup is a family $(\Pt)_{t\geq 0}$ of operators acting on bounded measurable functions on  $(E, \mathcal E)$, which preserve positivity and satisfy $\Pt(1)=1$ and the semigroup property $P_t\circ P_s= P_{t+s}$. We assume the existence of a reference measure $\mu$, in general chosen as to be an invariant measure for $(\Pt)_{t\geq0}$, and  we assume that $\Pt(f)$ converges to $f$ in $L^2(\mu)$ when $t\to 0$. In this context, the semigroup is entirely described by its infinitesimal generator $L,$ which is $\partial_t \Pt_{\vert t=0}$, so that $\Pt(f) $ solves the heat equation 
$\partial_t \Pt(f) = L\Pt(f)$. 

 We assume the existence of  a nice dense algebra $\mathcal A$ of functions on $E$ in the domain of $L$, and we define the {\it carr\'e du champ} operator $\Gamma$  by
$$
\Gamma(f,g)=\frac{1}{2}\PAR{L(fg)-fLg-gLf}, 
$$
for all functions $f$ and $g$ in $\mathcal A$. For simplicity we let $\Gamma(f) = \Gamma(f,f)$ denote the quadratic form. 

The operator $L$  is called  a {\it diffusion operator} if 
\begin{equation}
\label{eq-difff}
L\phi(f)=\phi'(f)Lf+\phi''(f)\Gamma(f)
\end{equation}
for  all $f\in\mathcal A$ and all $\phi\in\mathcal C^{\infty}(\dR).$ This means that the operator $L$ is a second order differential operator. This is equivalent to the fact  that  the carr\'e du champ operator $\Gamma$ satisfies the chain rule 
$$
\Gamma(\phi(f),g)=\phi'(f)\Gamma(f,g)
$$ 
for all functions $f$ and $g$,  or in other words that $\Gamma$ is a derivation in each variable.  

The classical example of diffusion operators is given in $\dR^n$ by 
$$
L f (x)= \sum_{i,j=1}^n D_{ij} (x) \frac{\partial^2 f}{\partial x_i \partial x_j} (x) - \sum_{i=1}^n a_i (x)  \frac{\partial f}{\partial x_i} (x)
$$
where  $D(x) = (D_{ij}(x))_{1 \leq i,j \leq n}$ is a symmetric $n \times n$ matrix, nonnegative in the sense of quadratic forms on $\dR^n$ and with smooth coefficients; also $a(x) = (a_i(x))_{1 \leq i \leq n}$ has smooth coefficients.  In this case one has :
$$
\Gamma (f) (x) = (D(x) \nabla f (x)) \cdot \nabla f(x). 
$$

A first fundamental instance is the heat semigroup $(\Pt)_{t \geq 0}$ on $\dR^n$, given by
\begin{equation}
\label{eq-heat}
\Pt f(x)=\int f(y)\exp\PAR{-\frac{|x-y|^2}{4t}}\frac{1}{(4\pi t)^{n/2}}dy,
\end{equation}
or equivalently 
$$
\Pt f(x)=E(f(x+\sqrt{2t}Y))
$$
where $Y$ is a random variable of law $\gamma$. Its generator is the Laplacian $\Delta$ on $\dR^n$.

Another instance of great interest is the Ornstein-Uhlenbeck semigroup $(\Nt)_{t\geq 0}$ on $\dR^n$, given by the Mehler formula  
\begin{equation}\label{eq-ExplicitOU}
\Nt(f)(x)= E(f(e^{-t} x+ \sqrt{1-e^{-2t}}Y)).
\end{equation}

Its generator is $L=\Delta-x\cdot\nabla.$ It is linked with the heat semigroup $(\Pt)_{t \geq 0}$ in the following way : Let us define the semigroup of dilations by
$$
\Tt f(x)=f(e^{t/2}x)
$$
for all $t\geq 0$ and functions $f$, so that
\begin{equation}\label{eq-comPT}
\TT{a}{\PT{b}{f}}(x)=E(f(e^{a/2}x+\sqrt{2b}Y)) = \PT{b e^{-a}}{\TT{a}{f}}(x).
\end{equation}
Then,  the heat and Ornstein-Uhlenbeck semigroups are related by
\begin{equation}\label{eq-PN}
\TT{a}{\PT{b}{f}}=\Nt(f)
\end{equation}
with $a=-2t$ and $b=(1-e^{-2t})/2.$

One can also deal on differentiable manifolds with the same kind of operators given in a local system of coordinates, and this is the case for example for the Laplace-Beltrami operator on a manifold.

%{\red fb: plus sur les vari\'et\'es ?}

\medskip

Assumptions on $L$ will be given in terms of the $\Gamma_2$ operator defined by
 $$
  \Gamma_2(f) = \frac{1}{2} \PAR{L \Gamma(f) - 2\Gamma (f, Lf)}
 $$
for all functions $f$ in $\mathcal A$. Although $\Gamma(f)$ is always non-negative (and this is in some way a characteristic property of Markov generators), it is not the case in general for the operator~$\Gamma_2$. 
% \begin{edefi}

If $n\geq 1$ and $\rho\in\dR$, we say that the semigroup $(\PT{t})_{t \geq 0}$ (or the infinitesimal generator $L$) satisfies the {\it curvature-dimension $CD(\rho,n)$  criterion} if
$$
\Gamma_2(f)\geq \rho\Gamma(f)+\frac{1}{n}\PAR{Lf}^2
$$
for all functions $f\in\A$.
%\end{edefi}
The $\Gamma_2$ operator and the associated  $CD(\rho,n)$ criterion have been introduced by the first author and M. Emery in \cite{bakryemery}. 

For instance, the Laplacian on $\dR^n$ satisfies the $CD(0,n)$ criterion and the operator $L = \Delta - a(x) \cdot \nabla$ satisfies the $CD(\rho, \infty)$ criterion if and only if $Ja(x) + (Ja(x))^* \geq 2 \rho Id$ as symmetric matrices on $\dR^n$, where $Ja$ is the Jacobian matrix of $a.$ In particular the Ornstein-Uhlenbeck semigroup satisfies the $CD(1, \infty)$ criterion. The Laplace-Beltrami operator on the $n$-dimensional unit sphere satisfies the $CD(n-1,n)$ criterion.

\medskip

One of the applications of the criterion is  to describe the local estimate of the logarithmic Sobolev constant of the semigroup $(\Pt)_{t\geq 0}$. Indeed, letting
$$
\ent{\Pt}{f}=\Pt (f\log f)-\Pt f\log\Pt f,
$$
the curvature-dimension $CD(\rho,\infty)$  criterion is equivalent to the local logarithmic Sobolev inequality 
\begin{equation}\label{logsob1}
%\label{localLSI}
\ent{\Pt}{f} \leq \frac{1 - e^{-2\rho t}}{2\rho} \, \Pt\PAR{\frac{\Gamma (f)}{f}}
\end{equation}
for all positive functions $f\in\mathcal A$ and $t \geq 0$. Moreover, for $\rho >0,$ it leads to the logarithmic Sobolev inequality 
$$
\ent{\mu}{f} \leq \frac{1}{2\rho} \, \mu\PAR{\frac{\Gamma (f)}{f}}
$$
for the ergodic measure $\mu,$ where $\ent{\mu}{f} = \mu(f \log f)-\mu(f)\log \mu(f).$

%{\red fb: ai ajout\'e cette phrase car elle est sous entendue qd apres le theoreme 3.1 on fait tendre t vers l'unfini pour recuperer Gross. Mais on peut l'enlever de ce paragraphe "local"}

On the other hand the curvature-dimension $CD(0,n)$  criterion is equivalent to the local logarithmic Sobolev inequality 
\begin{equation}\label{logsob2}
\ent{\Pt}{f} \leq tL\Pt f+\frac{n}{2} \Pt f \, \log\PAR{1-\frac{2t}{n}\frac{\Pt\PAR{fL(\log f)}}{\Pt f}}
\end{equation}
 for any positive functions $f\in \mathcal A$ and $t\geq 0$. Inequality~\eqref{logsob1} is proved in~\cite{bakrytata} and~\eqref{logsob2} in~\cite{bakry-ledoux06}.
 
 There is an analog of  (\ref{logsob2}) for the general condition $CD(\rho,n)$, but is has a far less pleasant form and we shall not use it in what follows.
 
 It is the purpose of the following two sections to show how, in turn, these local logarithmic Sobolev inequalities translate into hypercontractivity estimates on the semigroup. In the last section, we shall use a curvature-dimension criterion argument for  a L\'evy semigroup, which is a non-diffusive Markov semigroup. 

%%%%%%%%%%%%%%%%%%%%%%%%%%%%%%%%%%%%%
\section{Hypercontractivity for diffusions under $CD(\rho,\infty)$}
\label{sec-cdrho}
%%%%%%%%%%%%%%%%%%%%%%%%%%%%%%%%%%%%%

In this section, we exploit the local inequality (\ref{logsob1})  to obtain new equivalent forms of the $CD(\rho, \infty)$ criterion.  In particular, applied to the usual Euclidean heat kernel in $\dR^n$ (and not the Ornstein-Uhlenbeck one), they recover Nelson's theorem. Although the results given here are certainly not new, it is interesting to compare them with the results that we shall state in the next section under the $CD(0,n)$ inequality, in which Nelson's theorem shall be extended  using  a similar method.

The simplest way to give an interpretation of the $CD(0,\infty)$ inequality in terms of the semigroup $(\Pt)_{t \geq 0}$ is perhaps, for a given function $f$, to let
$$
\Phi(s)= \PT{s}(\PT{t-s}f)^2
$$ on $[0,t]$. Then,  by the very definitions of the operators $\Gamma$ and $\Gamma_2$,
$$
\Phi'(s)= 2\PT{s}(\Gamma(\PT{t-s}f)), ~\Phi''(s)= 4 \PT{s} (\Gamma_2(\PT{t-s}f)).
$$
Therefore the $CD(0,\infty)$ inequality just translates into the convexity of the function $\Phi$, which for example is equivalent  to 
$$\PT{s}(\PT{t-s}f)^2 \leq(1- \frac{s}{t})(\PT{t} f)^2 + \frac{s}{t} \PT{t}(f^2)$$
for all $s\in [0,t]$. One would have a similar inequality under $CD(\rho,\infty)$.  This inequality has the advantage of being quite stable under convergence of semigroups, since it applies to functions $f$ and not to their derivatives, like $\Gamma(f)$.

But, in view of the applications, this inequality is not very helpful, compared with inequality (\ref{eq-majrho}) below. In particular, for the Euclidean heat kernel or the Ornstein-Uhlenbeck semigroup, it would be quite difficult to compare it with the Nelson's hypercontractivity theorem. This is why we shall give another equivalent form of this criterium.

\begin{ethm}\label{theo-rho}
Let $\rho\in\dR$ and $(\PT{t})_{t \geq 0}$ be a diffusion Markov semigroup. 
Then the following assertions are equivalent :
\begin{enumerate}[(i)]
\item
 the semigroup $(\PT{t})_{t \geq 0}$ satisfies the ${CD}(\rho,\infty)$ criterion;
 \item 
 the semigroup $(\PT{t})_{t \geq 0}$ satisfies the local logarithmic Sobolev inequality :
\begin{equation}
\label{eq-localsob}
\ent{\Pt}{f} \leq \frac{1 - e^{-2\rho t}}{2\rho} \, \Pt\PAR{\frac{\Gamma (f)}{f}}
\end{equation}
 for all positive functions $f\in \mathcal A$ and $t\geq 0$; 
%  \item 
%For all  $0<q_2< q_1<1$, $0<s<t$, satisfying  
% \begin{equation}
%\label{eq-majrho2}
%{\PT{t}(f^{q_1})}^{1/q_1}\leq {\PT{s}((\PT{t-s}f)^{q_2})}^{1/q_2} 
%\end{equation}
%for all nonnegative function $f\in\mathcal A$,
 \item 
for all   $0<s \leq t$ and   $1<q_1 \leq q_2$ such that 
\begin{equation}
\label{eq-rela}
\frac{q_2-1}{q_1-1}=\frac{e^{2\rho t}-1}{e^{2\rho s}-1}
\end{equation}
one has
\begin{equation}
\label{eq-majrho}
{\PT{s}((\PT{t-s}f)^{q_2})}^{1/q_2}\leq {\PT{t}(f^{q_1})}^{1/q_1}
\end{equation}
for all positive functions $f\in\mathcal A$; 
 \item 
 for all $0<s \leq t$ and all $p$ and $q$ such that $0<q_2 \leq  q_1<1$ or  $q_2 \leq  q_1<0$, and~\eqref{eq-rela}, one has
\begin{equation}
\label{eq-majrho2}
 {\PT{t}(f^{q_1})}^{1/q_1}\leq{\PT{s}((\PT{t-s}f)^{q_2})}^{1/q_2}
\end{equation}
for all positive functions $f\in\mathcal A$.
 \end{enumerate}
 If $\rho=0$, then $\PAR{1 - e^{-2\rho t}}/\PAR{2\rho}$ will be replaced by $t$ and $\PAR{e^{2\rho t}-1}/\PAR{e^{2\rho s}-1}$ by ${t}/{s}$.  
\end{ethm}
\begin{eproof}
The equivalence between $(i)$ and $(ii)$ is proved in \cite{bakrytata}, one can also see~\cite[Chapter 5]{logsob}. 

Let us now assume $(ii)$ and prove $(iii)$ and $(iv)$. Let $f\in\mathcal A$ be a positive function and define the function $\psi$ on $[0,t]$ by 
$$
\psi(s)=\PT{s}{\PAR{\PT{t-s}{ \PAR{f}^{q}}}}^{1/q},
$$  
where  $q$ is a  function of $s$.  If  $g=\PT{t-s}{f}$ then 
$$
\psi'\psi^{q-1}\frac{q^2}{q'}=\ent{\PT{s}}{g^q}+\frac{q^2}{q'}({q-1})\PT{s}(g^{q-2}\Gamma(g)) \leq  q^2\PAR{\frac{1 - e^{-2\rho s}}{2\rho}+\frac{q-1}{q'}}\PT{s}\PAR{g^{q-2}\Gamma(g)}
$$
by the local logarithmic Sobolev inequality~\eqref{eq-localsob}. Let now $q$ satisfy $\frac{1 - e^{-2\rho s}}{2\rho}+\frac{q-1}{q'}=0$, that is,
$$
\frac{q(s)-1}{q(t)-1}=\frac{e^{2\rho t}-1}{e^{2\rho s}-1}
$$
for all $s\in(0,t]$. Depending on whether $q$ is decreasing or increasing, or equivalently on whether $q>1$ or $q<1$, the function $\psi$ is increasing or decreasing. The two cases give the two inequalities~\eqref{eq-majrho} and~\eqref{eq-majrho2}. 
\smallskip

Then we assume $(iii)$. We apply a first order Taylor expansion for $s$ going to $t$ and $q_2$ going to $q_1$: we let  $q_1=2$, $q_2=2(1+ \varepsilon)$, and $s=t(1-\alpha \varepsilon)+o( \varepsilon)$ for $ \varepsilon>0$ and $\alpha=2(1-e^{-2\rho t})/(2\rho t)$. This relation implies that the condition~\eqref{eq-rela} is satisfied as $ \varepsilon$ goes to $0$.  The Taylor expansion of the inequality~\eqref{eq-majrho} gives 
$$
\ent{\Pt}{f^2}\leq 4\frac{1-e^{-2\rho t}}{2\rho}\Pt(\Gamma(f)),
$$
which is inequality~\eqref{eq-localsob} for $f^2$ instead of $f.$ The same argument ensures that $(iv)$ implies $(ii)$. 
\end{eproof}

\begin{erem}
Inequalities~\eqref{eq-majrho} and~\eqref{eq-majrho2} give criteria equivalent to the $CD(\rho,\infty)$ criterion without resorting to any derivation;  in particular they apply to all bounded measurable functions (whence their stability under convergence of semigroups). 
\end{erem}

\begin{erem}\label{rem-gross}
Under the condition $CD(\rho,\infty)$ with $\rho >0,$ Theorem~\ref{theo-rho} leads to the hypercontractive bound given by the Gross Theorem. Let indeed a semigroup $(\Pt)_{t\geq0}$ satisfy the $CD(\rho,\infty)$ criterion and  be ergodic in $L^p(d\mu)$ where $\mu$ is an invariant probability measure. If $s=t-u$ with $u \geq  0$  then inequality~\eqref{eq-majrho} gives 
$$
{\PT{t-u}((\PT{u}f)^{q_2})}^{1/q_2}\leq {\PT{t}(f^{q_1})}^{1/q_1},
$$  
where $1<q_1 \leq q_2$ satisfy~\eqref{eq-rela}. Now, letting $t$ go to infinity, we recover the classical hypercontractive bound
$$
\NRM{\PT{u}f}_{L^{q_2}(d\mu)}\leq \NRM{f}_{L^{q_1}(d\mu)},
$$  
where $1<q_1 \leq q_2$ satisfy $q_2-1=e^{2\rho u}(q_1-1)$. Analogously from~\eqref{eq-majrho2}, we recover the reverse hypercontractive bound
\begin{equation}
\label{eq-derder}
\NRM{f}_{L^{q_1}(d\mu)}\leq \NRM{\PT{u}f}_{L^{q_2}(d\mu)},
\end{equation}
where $0<q_2 \leq q_1<1$ or $q_2 \leq q_1<0$ satisfy $q_2-1=e^{2\rho u}(q_1-1)$.
\end{erem}

%{\red Il faudrait faire un paragraphe, la, ou enoncer qqch. proposition:}

\begin{erem}
\label{rem-ou}
According to the above remark, Theorem~\ref{theo-rho} implies Nelson's Theorem when applied to the Ornstein-Uhlenbeck semigroup which satisfies the $CD(1, \infty)$ inequality.  But there is another way of recovering this theorem, using the heat equation in the Euclidean space.  This way may seem spurious at first glance, but we shall see in the next section that it carries much more information.    

As pointed out in Section~\ref{sec-diff}, the Laplacian operator on $\dR^n$ satisfies the $CD(0, \infty)$ criterion. Hence we may apply the bound~\eqref{eq-majrho}  to the heat semigroup $(\Pt)_{t\geq 0}$, at $x=0$, in the form
\begin{equation}
\label{eq-ou2}
\PT{s}((\PT{t-s}f)^{q_2})^{1/{q_2}} (0) \leq \Pt (f^{q_1})^{1/{q_1}} (0)
\end{equation}
where $q_2=1+\frac{t}{s}(q_1-1)$, that is,
$$
\NRM{\TT{\log(2s)}{\PT{t-s}{f}}}_{L^{q_2}(d\gamma)} \leq  \NRM{\TT{\log(2t)} f}_{L^{q_1}(d\gamma)}
$$
by~\eqref{eq-heat}, or
$$
\NRM{\TT{\log(2s)}{\PT{t-s}\TT{-\log(2t)}{g}}}_{L^{q_2}(d\gamma)} \leq  \NRM{g}_{L^{q_1}(d\gamma)}
$$
for all positive $g.$ But, by~\eqref{eq-comPT} and~\eqref{eq-PN},
$$
\TT{\log(2s)}{\PT{t-s}\TT{-\log(2t)}{g}} = \TT{\log(s/t)}{\PT{(t-s)/(2t)}{g}} = \Nt g
$$
for $s = t e^{-{2t}}$, which leads to the bound
$$
\NRM{\Nt g}_{L^{q_2}(d\gamma)}\leq \NRM{g}_{L^{q_1}(d\gamma)},
$$
given by Nelson's Theorem, with $q_2-1= e^{2t}(q_1-1)$.
\end{erem}

%{\red fb: la demo est maintenant plus courte car tout est dans la section 2}

\begin{erem} The relations  through dilations between the heat kernel and the Ornstein-Uhlenbeck kernel on $\mathbb R^n$ just reflect the commutation properties of their  generators. Indeed, if $D$ is the operator defined by
$$D(f)= \frac{1}{2} \sum_{i=1}^n  x_i \frac{\partial f}{\partial x_i},
$$ then
$$T_t f= \exp(tD),$$ while
$$[\Delta, D]= \Delta ~\hbox{ and } L= \Delta-2D.$$
A similar relation holds in the hypoelliptic setting for stratified nilpotent groups, for example the Heisenberg  group with the Kohn Laplacian, but in these situations there is no $CD(\rho, \infty)$ inequality, and the $\Gamma_2$ calculus is much more delicate to use, see for example \cite{Bakry-Baudoin-Bonnefont-Qian}. It is not clear that one could obtain similar results in this context.

\end{erem}

%%%%%%%%%%%%%%%%%%%%%%%%%%%%%%%
%%%%%%%%%%%%%%%%%%%%%%%%%%%%%%%
\section{Hypercontractivity for diffusions under $CD(0,n)$}
\label{sec-cdn}
%%%%%%%%%%%%%%%%%%%%%%%%%%%%%%%
%%%%%%%%%%%%%%%%%%%%%%%%%%%%%%%

In this section, we extend the equivalences of the previous one to the case of the $CD(0,n)$ inequality. 

\subsection{A general dimensional hypercontractivity condition}

\begin{ethm}\label{theo-cdn}
Let $n\geq 1$ and $L$ be a diffusion Markov semigroup. 
Then the following assertions are equivalent :
\begin{enumerate}[(i)]
\item
 the semigroup $(\PT{t})_{t \geq 0}$ satisfies the ${CD}(0,n)$ criterion;
 \item 
 the semigroup $(\PT{t})_{t \geq 0}$ satisfies the local logarithmic Sobolev inequality :
\begin{equation}
\label{eq-logsobn}
\ent{\Pt}{f} \leq tL\Pt f+\frac{n}{2} \Pt f  \log\PAR{1-\frac{2t}{n}\frac{\Pt\PAR{fL(\log f)}}{\Pt f}}
\end{equation}
 for all positive functions $f\in \mathcal A$ and $t\geq 0$; 
 \item 
for all $1<q_1< q_2$, $u_2,u_1\geq 0$ and  $t-s=u_1q_1-u_2q_2$, 
\begin{equation}
\label{eq-maj}
{\PT{u_2}((\PT{t-s}f)^{q_2})}^{1/q_2}\leq {\PT{u_1}(f^{q_1})}^{1/q_1}{M}^{n/2}
\end{equation}
for all positive functions $f \in \mathcal A$ where 
$$
M=\PAR{\frac{q_1-1}{u_2}}^{1-1/q_1}\PAR{\frac{q_2-1}{u_1}}^{1/q_2-1}
\PAR{\frac{u_1q_1-u_2q_2}{q_2-q_1}}^{1/q_2-1/q_1};
$$
 \item 
For all  $0<q_2< q_1<1$ or  $q_2< q_1<0$ and  $u_2,u_1\geq 0$ satisfying  $t-s=u_1q_1-u_2q_2$, 
 \begin{equation}
\label{eq-maj2}
{\PT{u_1}(f^{q_1})}^{1/q_1}\leq {\PT{u_2}((\PT{t-s}f)^{q_2})}^{1/q_2} {N}^{n/2}
\end{equation}
for all positive functions $f\in\mathcal A$, where 
$$
N=\PAR{\frac{1-q_1}{u_2}}^{1/q_1-1}\PAR{\frac{1-q_2}{u_1}}^{1-1/q_2}
\PAR{\frac{u_1q_1-u_2q_2}{q_1-q_2}}^{1/q_1-1/q_2}.
$$
 \end{enumerate}
  \end{ethm}

\begin{erem}
Inequalities~\eqref{eq-maj} and~\eqref{eq-maj2} are optimal in the sense that if $L$ is the Laplacian in $\dR^n$ then these inequalities are equalities for square-exponential (or Gaussian) functions.

Indeed, the Laplacian operator is a diffusion in $\dR^n$ satisfying the $CD(0,n)$ criterion. The associated heat semigroup $(\Pt)_{t \geq 0}$ is given by~\eqref{eq-heat}. In particular, if $f(x)=e^{a|x|^2}$ with $a\in\dR$, then
$$
\Pt f(x)=\frac{1}{\PAR{1-4ta}^{n/2}}\exp\PAR{\frac{a|x|^2}{1-4ta}}
$$
 for all $x\in\dR^n$ and $t$ such that $1-4ta>0$ so that
$$
\PAR{\frac{{\PT{u_2}((\PT{t-s}f)^{q_2})}^{1/q_2}(x)}{ {\PT{u_1}(f^{q_1})}^{1/q_1}(x)}}^{2/n}=\frac{\PAR{1-4q_1u_1a}^{1/q_1}}{{1-4a(t-s)}}\PAR{\frac{1-4a(t-s)}{1-4q_1u_1a}}^{1/q_2}
$$
 for all $x\in\dR^n$ since $t-s=u_1q_1-u_2q_2$. It follows that  the inequality~\eqref{eq-maj} becomes an equality for all $x$, provided we choose
$$
a=\frac{t-s+u_2-u_1}{4 u_1(t-s)(q_1-1)}.
$$

%{\red fb: il manquait le $u_1$ au denominateur, a verifier  IG: OK}

E.H.  Lieb in~\cite{lieb} proves that operator bounds on Gaussian kernels as operators from $L^p(\mathbb R^n)$ to $L^q(\mathbb R^n)$ have only Gaussian maximizers. It seems to us that it should also be the case for~\eqref{eq-maj} in the setting of the Laplacian in $\dR^n$. But the inequality~\eqref{eq-maj} is more complicated than an operator bound and we do not know how to prove it.

\end{erem}

\textbf{\textit{Proof of Theorem~\ref{theo-cdn}}} 
\noindent \proofbegin\,
For future use we note that the bound \eqref{eq-logsobn} in $(ii)$ is equivalent to the family of bounds
\begin{equation}\label{eq-logsobn'}
 \ent{\Pt}{f} \leq tL\Pt f+\frac{n}{2} (\lambda - 1 - \log \lambda) \Pt f  - t \, \lambda \, \Pt\PAR{fL(\log f)}
 \end{equation}
 for $\lambda >0,$ since
 $$
\log (1+x) = \sup_{\lambda >0} [ \lambda x +\lambda - 1 - \log \lambda ]
$$
 by concavity of the logarithm.

The equivalence between $(i)$ and $(ii)$ is proven  by the first author and M. Ledoux in \cite[Theorem 1]{bakry-ledoux06}. 

Let us now assume $(ii)$ and prove $(iii)$ and $(iv)$. Let $f\in\mathcal A$ be a positive function and define the function $\psi$ of $s\in[0,t]$ by 
$$
\psi(s)=\PT{u}{\PAR{\PT{t-s}{ \PAR{f}^{q}}}}^{1/q},
$$  
where  $q$ and  $u$ are functions of $s$. The function $u$ has to be positive. Then 
\begin{equation}
\label{eq-in1}
\psi'\psi^{q-1}\frac{q^2}{q'}=\ent{\PT{u}}{g^q}+\frac{q^2}{q'}(u'-1)\PT{u}(g^{q-1}Lg)+\frac{q^2}{q'}u'({q-1})\PT{u}(g^{q-2}\Gamma(g))
\end{equation} 
where $g=\PT{t-s}{f}.$ By~\eqref{eq-logsobn'}, the logarithmic Sobolev inequality~\eqref{eq-logsobn} applied to $g^q$ implies 
$$
\ent{\PT{u}}{g^q}\leq qu(1-\lambda)\PT{u}(g^{q-1}Lg)+qu(q-1+\lambda)\PT{u}(g^{q-2}\Gamma(g))+\frac{n}{2}A_\lambda{\PT{u}}\PAR{g^q}.
$$ 
for all $\lambda>0$, where $A_{\lambda} = \lambda - 1 - \log \lambda.$ Then~\eqref{eq-in1} becomes 
\begin{multline*}
\psi'\psi^{q-1}\frac{q^2}{q'}\leq \SBRA{qu(1-\lambda)+\frac{q^2}{q'}(u'-1)}\PT{u}(g^{q-1}Lg)\\
+\SBRA{qu(q-1+\lambda)+\frac{q^2}{q'}u'({q-1})}\PT{u}(g^{q-2}\Gamma(g))+\frac{n}{2}A_\lambda\psi^q.
\end{multline*}
Let now $q = q(s), u = u(s)$ and $\lambda = \lambda(s)$ solve
\begin{equation}
\label{eq-sys}
\left\{
\begin{array}{l}
\disp\frac{q}{q'}(u'-1)+u(1-\lambda)=0\\
\disp u(q-1+\lambda)+\frac{q}{q'}u'({q-1})=0
\end{array}
\right.
\end{equation}
Then one has 
$$
\psi'\psi^{-1}\frac{q^2}{q'}\leq \frac{n}{2}A_\lambda.
$$

\noindent
{\it First case : $q$ is a decreasing function.} 
Then, by integration,
$$
\frac{\psi(t)}{\psi(s)}\geq \exp\PAR{\frac{n}{2}\int_s^t\frac{q'(r)}{q(r)^2}A_{\lambda(r)} \, dr }:= \exp \PAR{\frac{n}{2}M_{s,t}},
$$
that is,
$$
{\PT{u(s)}(\PT{t-s} (f)^{q(s)})}^{1/q(s)}\leq \exp\PAR{-\frac{n}{2}M_{s,t}}{\PT{u(t)}(f^{q(t)})}^{1/q(t)}
$$
The two relations between $\lambda$, $u$ and $q$ give  
$$
(uq)' (r)=1
$$
for $s \leq r \leq t.$ In other words there exists a constant $C$ such that  
$$
u(r)=\frac{r+C}{q(r)}, 
$$
and therefore, by the first equation in~\eqref{eq-sys}, 
$$
\lambda=\frac{1-q}{q'u}=\frac{q(1-q)}{q'(r+C)}.
$$
In particular, if $q$ is decreasing, then  $q>1$ (we always have $u(r)={q(r)}/{(r+C)}\geq 0$). 

If  $t>s>0$ are fixed, the quantity $M_{s,t}$ depends only on the function $q$ on $[s,t]$ and the constant $C$. Let now the parameters $t$, $s$, $u(t)=u_1$,  $q(t)=q_1$, $q(s)=q_2$ and $q(t)=q_1>1$ be fixed. Then the quantity 
$$
M_{s,t} = \int_s^t \frac{1-q}{q(r+C)} - \frac{q'}{q^2} - \frac{q'}{q^2} \log \frac{q(1-q)}{q(r+C)} (r) \, dr
$$ 
is made maximal by a map $q$ such that
$$
\Big(\frac{q^2}{(r+C)^2 q')} \Big)' = 0
$$
on $[s,t]$, that is given on $[s,t]$ by 
$$
 q(r)=\frac{r+C}{\alpha r+\beta},
 $$
where $\beta$ is a constant and 
 $$
 \alpha=\frac{u_2-u_1}{s-t}=\frac{u_2-u_1}{u_2q_2-u_1q_1}.
 $$
Then
$$
M_{s,t}= \int_s^t\frac{q'(r)}{q(r)^2}A_{\lambda(r)} \, dr = \int_{q_1}^{q_2}\frac{1}{q^2}\PAR{\frac{1-q}{1-\alpha q}-1-\log\frac{1-q}{1-\alpha q}}dq
$$ 
 since  $ \lambda=\PAR{1-q}/\PAR{1-\alpha q}$ for this choice of $q.$  Moreover an antiderivative of
 $$
\frac{1}{q^2}\PAR{\frac{1-q}{1-\alpha q}-1-\log\frac{1-q}{1-\alpha q}}
$$ 
is given by $-\frac{q-1}{q}\log\frac{q-1}{\alpha q-1}$.  Hence  
$$
\exp\PAR{-\frac{n}{2} \int_{s}^t\frac{q'}{q^2}A_\lambda dr}=\SBRA{\PAR{\frac{q_1-1}{u_2}}^{1-1/q_1}\PAR{\frac{q_2-1}{u_1}}^{1/q_2-1}
\PAR{\frac{u_1q_1-u_2q_2}{q_2-q_1}}^{1/q_2-1/q_1}}^{n/2},
$$
which proves the assertion $(iii)$. 

\noindent
{\it Second case : $q$ is an increasing function.}
The computation is the same as in the previous case.  In this case one has to assume that $q\in(0,1)$ or $q<0$. One obtains $(iv).$

\medskip

Let us now prove that $(iii)$ implies $(ii)$.  As usual we perform a first order Taylor expansion in inequality~\eqref{eq-maj}. We let $f$ be a positive function in the algebra $\mathcal A$ and let  $q_1 = 2, q_2=2(1+ \varepsilon), u_1 = t$ and $u_2=t(1-a \varepsilon)$ for some $t, \varepsilon,a>0$, so that $s = t(1 + 2(1-a)  \varepsilon) + o( \varepsilon).$ A Taylor expansion for $ \varepsilon$ going to $0$ leads to
$$
M^{n}=1+ \frac{ \varepsilon n}{2}(a - 2 - \log (a-1)) + o( \varepsilon)
$$
and then to 
$$
\ent{\Pt}{f^2} + t(a-2) \Pt L(f^2) - 4t(a-1) \Pt (\Gamma(f)) \leq \frac{n}{2} (a-2 - \log(a-1)) \Pt (f^2)
$$
that is,
$$
\ent{\PT{t}}{g}\leq t(1-\lambda)\PT{t}(Lg)+\frac{n}{2}(\lambda-1-\log \lambda)\PT{t}(g)+t\lambda\PT{t}\PAR{\frac{\Gamma(g)}{g}}.
$$
with $g=f^2$ and $\lambda = a - 1.$ Since $a >1$ is arbitrary, so is $\lambda >0$, which leads to~\eqref{eq-logsobn} by~\eqref{eq-logsobn'}. The same computation can be done starting from $(iv)$ instead of $(iii)$, which concludes the proof of the theorem. 
\proofend

\bigskip

For all $\lambda>0$,  then $\lambda-1-\log\lambda\geq 0$ and it is equal to 0 if and only if $\lambda=1$. As we can see in the proof it implies that we always have $M\geq 1$,  and $M=1$ if and only if $u' \equiv 1$. With this particular choice of $u_1 = t$ and $u_2 = s$ we recover the result of Theorem~\ref{theo-rho} with $\rho=0$, for which the hypercontractive bound does not depend on the dimension.

\bigskip

\subsection{The case of the Ornstein-Uhlenbeck semigroup}\label{sec-appliOU}

We have seen in Remark~\ref{rem-ou} how to recover the Nelson Theorem only using a $CD(0, \infty)$ inequality on the Euclidean heat semigroup. In fact the semigroup satisfies the stronger $CD(0, n)$ condition, whence the inequality~\eqref{eq-maj}.  Following the argument in Remark~\ref{rem-ou}, the Ornstein-Uhlenbeck can be seen as a combination of the dilation and heat semigroups, and we obtain
\begin{equation}
\label{eq-nash1}
\NRM{\NT{t}{f}}_{L^{q_2}(d\gamma)}\leq M^{n/2}\NRM{\TT{-a}{f}}_{L^{q_1}(d\gamma)}
\end{equation}
for the Ornstein-Uhlenbeck semigroup $(\Nt)$, where
$$
M=\PAR{\frac{q_1-1}{e^{-2t}}}^{1-1/q_1}\PAR{\frac{q_2-1}{e^{-a}}}^{1/q_2-1}\PAR{\frac{1-e^{-2t}}{q_2-q_1}}^{1/q_2-1/q_1}
$$
for all $t\geq 0, a \in \dR$ and $1< q_1<q_2$ satisfying $q_2-1=e^{2t}(q_1e^{-a}-1).$

When $a=0$, then $M=1$ and inequality~\eqref{eq-nash1} is simply the Nelson  classical hypercontractivity estimate~\eqref{eq-hyper1} of the Ornstein-Uhlenbeck semigroup.  We cannot obtain a better inequality in terms of a classical hypercontractivity property : if $q_1$ and $t$ are given, there is no larger $q_2$ such that~\eqref{eq-hyper1} holds even with an additional constant.

The main difference is that the limit $q_1$ going to $1$ is not informative in Nelson's Theorem. On the contrary, in inequality~\eqref{eq-nash1}, and for $q_2 = 2$, this limit  leads to the following inequality for Gaussian kernels
\begin{equation}\label{Nash1OU}
\NRM{\NT{t}{f}}_{L^{2}(d\gamma)}\leq \PAR{\frac{e^{-2t}+1}{1-e^{-2t}}}^{n/4}\NRM{\TT{\ln\PAR {e^{-2t}+1}}{f}}_{L^{1}(d\gamma)}
\end{equation}
for $t>0$.  The coefficient in~\eqref{Nash1OU} behaves like $t^{-n/4}$ for small $t$, which can be related to the hypercontractive bound
$$
\NRM{\PT{t}{f}}_{L^{2}(dx)}\leq \frac{1}{(8 \pi t)^{n/4}} \NRM{f}_{L^{1}(dx)}
$$
for the heat semigroup and with respect to the Lebesgue measure.

Although the Ornstein-Uhlenbeck semigroup does not satisfy any $CD(\rho, n)$ inequality, its relation through dilations with the Euclidean heat semigroup, which satisfies $CD(0,n)$, carries dimensional properties of this semigroup. In particular, inequality (\ref{Nash1OU}) implies the finitness of the trace of $\Nt$  (and in fact its exact value), which depends on the dimension and cannot be deduced from the sole $CD(1, \infty)$ property of its generator.

Indeed, by a standard change of variables, we have
$$
\PAR{\frac{e^{-2t}+1}{1-e^{-2t}}}^{n/4}\NRM{\TT{\ln\PAR {e^{-2t}+1}}{f}}_{L^{1}(d\gamma)} =\int |f(y)|V_t(y) \gamma(dy),
$$ 
where  
\begin{equation}\label{Nash4OU}
V_t(y)=(1-e^{-4t})^{-n/4}\exp \Big(\frac{|y|^2}{2}\frac{1}{1+ e^{2t}} \Big).
\end{equation}
The inequality (\ref{Nash1OU}) now reads
\begin{equation}\label{Nash2OU}
\NRM{\NT{t}{f}}_{L^{2}(d\gamma)}\leq \int |f|V_t d\gamma.
\end{equation}

But, if a  Markov operator $K$, symmetric with respect to a measure $\mu,$ satisfies 
\begin{equation} \label{Nash3OU}
\|Kf\|_{L^2(d\mu)} \leq \| fV\|_{L^1(d\mu)}
\end{equation}
for a positive  function $V$, then the Markov kernel of the operator $K^2= K\circ K$ has a density  $k^2(x,y)$ with respect to  $\mu$ bounded from above by $V(x)V(y)$.  Let us indeed consider the operator 
$$
K_1(f) = \frac{1}{V}K(fV),
$$
 which is symmetric with respect to the measure $d\nu = V^2 d\mu$. Equation~\eqref{Nash3OU} reads
$\|K_1 f\|_2 \leq \|f\|_1$, the norms being now computed with respect to the measure $\nu$.
 A standard duality argument implies that $\|K_1 f\|_\infty \leq \|f\|_2$, so that $\|K_1^2 f\|_\infty  \leq \|f\|_1$;  in turn this implies that $K_1^2$ has a density with respect to $\nu$ which is bounded from above by $1$. This is equivalent to the fact that $K^2$ has a density  with respect to $\mu$ bounded from above by $V(x)V(y)$.
 
Back to the semigroup $(\Nt)_{t \geq 0}$, it follows from inequality~\eqref{Nash2OU}  that 
 $\NT{2t}$ has a density with respect to $\gamma$ bounded from above by $V_t(x)V_t(y)$. But for a symmetric Markov semigroup, such a bound is equivalent to the fact that the kernel $n_{2t}(x,y)$ is such that 
 $$
 n_{2t}(x,x) \leq V_t(x)^2.
 $$
But in fact $n_{2t}(x,x) = V_t(x)^2$ as ensured by comparing \eqref{Nash4OU} and the explicit value of $n_t$ given, for instance, by the Mehler formula~\eqref{eq-ExplicitOU}. Therefore, the obtained bound on $n_t(x,x)$ is sharp and so is the inequality~\eqref{Nash2OU}.  

It also implies the sharp bound $(1-e^{-t})^n$ on the trace 
 $$ 
 \int n_t(x,x) d\gamma (x) = \int V_{t/2}^2 d\gamma(x)
 $$ 
 of $\Nt$, which can be explicitly computed from~\eqref{Nash4OU}. Note that  it is dimension dependent, unlike the hypercontractivity bounds or the logarithmic Sobolev constant of the ergodic measure $\gamma.$

\medskip 
 
This illustrates the optimality of our method. 
Inequality (\ref{Nash2OU})  is in fact related to a weighted Nash inequality relative to the Ornstein-Uhlenbeck semigroup. This method of weighted Nash inequalities may be pushed forward to get information on the trace of various semigroups. This point of view shall be developed in the forthcoming article~\cite{bbg2}.

%%%%%%%%%%%%%%%%%%%%%%%%%%%%%%%%%%%%%%%%%%%%%%%%%%%%%%
\section{Hamilton-Jacobi equations and transportation inequalities} 
\label{sec-transport}
%%%%%%%%%%%%%%%%%%%%%%%%%%%%%%%%%%%%%%%%%%%%%%%%%%%%%%
\subsection{Hypercontractivity for Hamilton-Jacobi equations}
\label{sec-hj}

Solutions to Hamilton-Jacobi equations can be seen as limits of linear diffusion semigroups. Let indeed $L$ be an infinitesimal diffusion generator associated to  a Markov semigroup~$\PAR{\PT{t}}_{t\geq 0}$. 

Given $ \varepsilon>0$, let 
$$ 
u^\varepsilon  = \PT{\bf \varepsilon t} \big ( {\rm e}^{- f /2\varepsilon} \big )
$$
be the solution to ${\partial_t u^\varepsilon } = \varepsilon\, \GI u^\varepsilon$, with initial value $ {\exp}\PAR{- f /2\varepsilon }$. Then, since  $L$ is a diffusion generator, the map
$$
v^\varepsilon = - 2\varepsilon \log \PT{\bf\varepsilon t} \big ( {\rm e}^{- f
  /2\varepsilon } \big )
  $$
  is a solution of the  initial value partial differential equation
\begin{equation*}
\left\{ 
\begin{array}{rl}
\displaystyle
\frac{\partial v^\varepsilon }{ \partial t}
    + \frac{1}{ 2} \Gamma(v^\varepsilon)
  - \varepsilon\, \GI v^\varepsilon &= 0 \quad {\rm on} \,\, E \times
(0,\infty ), \\
\displaystyle
 v^\varepsilon  & = f \quad \! {\rm on} \,\, E \times \{t=0\},
\end{array}
\right. 
\end{equation*}
where $\Gamma$ is the {\it carr\'e du champ} associated to $L$. Now, as $\varepsilon \to 0$, it is expected that $v^\varepsilon $ converges to the solution $v$ of the Hamilton-Jacobi equation 
\begin{equation}
\label{eq-hj1}
\left\{ 
\begin{array}{rl}
\displaystyle
\frac{\partial v}{ \partial t}
    + \frac{1}{ 2} \Gamma(v) &
  = 0 \quad {\rm on} \,\, E \times
(0,\infty ), \\
\displaystyle
 v  & = f \quad \! {\rm on} \,\, E \times \{t=0\}.
\end{array}
\right. 
\end{equation}

%Given a smooth function $f$, and $ \varepsilon>0$, let $v^\varepsilon =
%v^\varepsilon (x,t)$ be the solution of the initial value partial differential equation
%\begin{equation*}
%\left\{ 
%\begin{array}{rl}
%\displaystyle
%\frac{\partial v^\varepsilon }{ \partial t}
%    + \frac{1}{ 2} \Gamma(v^\varepsilon)
%  - \varepsilon \GI v^\varepsilon &= 0 \quad {\rm on} \,\, E \times
%(0,\infty ), \\
%\displaystyle
% v^\varepsilon  & = f \quad \! {\rm on} \,\, E \times \{t=0\},
%\end{array}
%\right. 
%\end{equation*}
%where $\Gamma$ is the {\it carr\'e du champ} associated to $L$.
%As $\varepsilon \to 0$, it is expected that $v^\varepsilon $ approaches in a
%reasonable sense the solution of the Hamilton-Jacobi equation 
%\begin{equation}
%\label{eq-hj1}
%\left\{ 
%\begin{array}{rl}
%\displaystyle
%\frac{\partial v}{ \partial t}
%    + \frac{1}{ 2} \Gamma(v) &
%  = 0 \quad {\rm on} \,\, E \times
%(0,\infty ), \\
%\displaystyle
% v  & = f \quad \! {\rm on} \,\, E \times \{t=0\}.
%\end{array}
%\right. 
%\end{equation}
%On the other hand, since $L$ is a diffusion generator, the function $u^\varepsilon = {\rm e}^{-v^\varepsilon /2\varepsilon }$ satisfies the equation ${\partial u^\varepsilon \over \partial t} = \varepsilon \GI u^\varepsilon $ (with initial value $ {\rm e}^{- f /2\varepsilon }$).
%Therefore,
%$$ 
%u^\varepsilon  = \PT{\bf \varepsilon t} \big ( {\rm e}^{- f /2\varepsilon }
%\big ),
%$$
%which implies that the solution of~\eqref{eq-hj1}, at least formally,  is given by  
%$$
%v =\lim_{ \varepsilon\rightarrow0} - 2\varepsilon \log \PT{\bf\varepsilon t} \big ( {\rm e}^{- f
%  /2\varepsilon } \big ).
%$$ 
The solution at time $t$ of~\eqref{eq-hj1} will be denoted $\Qt f,$ which defines  a non-linear  semigroup $(\Qt)_{t\geq 0}$.  One can  see \cite{barles,evans} for a review on Hamilton-Jacobi equations.
\medskip

The classical example is given by $\Gamma(f)=|\nabla f|^2$ on a complete Riemannian manifold $(M,d)$, for which the Hamilton-Jacobi semigroup is explicitly given by the Hopf-Lax formula
\begin{equation}
\label{eq-hjma}
\Qt f (x)=\inf_y\BRA{f(y)+\frac{1}{2t}{d(x,y)}^2}, 
\end{equation}
when $f$ is a Lipschitz  function on $M$.  This fundamental example will be treated in section~\ref{sec-transine}.

\medskip

The Hamilton-Jacobi semigroup is strongly related to logarithmic Sobolev inequalities in the following way :  for a probability measure $\mu$ and a constant $C$, the logarithmic Sobolev inequality 
\begin{equation*}
%\label{eq-islf}
\ent{\mu}{f^2}\leq {2C}\int \Gamma( f)d\mu
\end{equation*}
for all functions $f$ is equivalent to the hypercontractivity of the Hamilton-Jacobi equation
\begin{equation}
\label{eq-bgl}
\NRM{e^{\Qt f}}_{L^{a+t/C}(d\mu)}\leq\NRM{e^{f}}_{L^{a}(d\mu)}
\end{equation}
for all $a>0$ and all functions $f$. This result has been proved in~\cite{bgl} and extended in~\cite{gentil02,gentil03}.  
%These equations has been studied in term of hypercontractivity properties and logarithmic Sobolev 
%inequality in~\cite{gentil02,gentil03}.

\medskip

We now turn from the invariant measure level of~\eqref{eq-bgl} to the local level,  to investigate how  a local hypercontractive bound on Hamilton-Jacobi equations can be obtained by using a curvature-dimension criterion $CD(\rho,n)$.

The first result concerns the  $CD(\rho,\infty)$  criterion:

\begin{ethm}\label{theo-hjrho}
Let $\rho\in\dR$ and $(\PT{t})_{t \geq 0}$ be a diffusion Markov semigroup. 
Then the following assertions are equivalent :
\begin{enumerate}[(i)]
\item
 the semigroup $(\PT{t})_{t \geq 0}$ satisfies the ${CD}(\rho,\infty)$ criterion;
 \item 
 for all $q_1>q_2> 0$ and  $u, t>0$  satisfying 
 \begin{equation}
 \label{eq-hjcondi}
 q_1=q_2+\frac{\rho t}{1-e^{-2\rho u}},
 \end{equation}
 then 
\begin{equation}
\label{eq-hjrho}
\PT{u}\PAR{e^{q_1\Qt f}}^{1/q_1}\leq \PT{u}\PAR{e^{q_2f}}^{1/q_2}
\end{equation}
for all functions $f$;
 \end{enumerate}
  If $\rho=0$ then~\eqref{eq-hjcondi} will be replaced by  $ q_1=q_2+t/(2u)$.
  \end{ethm}

   In the case when $\rho >0$ we note that~\eqref{eq-hjrho} leads to the bound~\eqref{eq-bgl} for the ergodic measure $\mu$ of the semigroup, by letting $u$ go to infinity.
   \medskip
  
   \begin{eproof}
    We just give the sketch of the proof since it follows the argument of the next theorem, which is given in greater detail.
 
 We first assume $(i)$ and prove $(ii)$. Given $u>0$ fixed we let $H(s) =  \PT{u}\PAR{e^{q(s) \QT{s} f}}^{1/q(s)}$ for $s$ on $[0,t].$ We differentiate with respect to $s$, use the local bound~\eqref{eq-localsob} at time $u$ and obtain the inequality $H'(s) \leq 0$ on $[0,t]$ provided 
 $$
 q(t) = q(s) + \frac{\rho}{1-e^{-2\rho u}}(t-s).
 $$
 This proves $(ii).$
 
 Then we assume $(ii)$. As $ \varepsilon$ goes to 0, a first order Taylor expansion of~\eqref{eq-hjrho}  with  $q_2=1$ and $t = \varepsilon$ gives the local logarithmic Sobolev inequality~\eqref{eq-localsob}, which implies $(i)$. 
 \end{eproof}

The second result concerns the $CD(0,n)$  criterion :

\begin{ethm}\label{theo-hjcdn}
Let $n\geq 1$ and $(\PT{t})_{t \geq 0}$ be a diffusion Markov semigroup. 
Then the following assertions are equivalent :
\begin{enumerate}[(i)]
\item
 the semigroup $(\PT{t})_{t \geq 0}$ satisfies the ${CD}(0,n)$ criterion;
 \item 
 for all $q_1>q_2> 0$, $u_2,u_1 > 0$ and  $t >0$  satisfying $t=2(u_1q_1-u_2q_2)$, 
\begin{equation}
\label{eq-hj2}
\PT{u_1}\PAR{e^{q_1\Qt f}}^{1/q_1}\leq \PT{u_2}\PAR{e^{q_2f}}^{1/q_2}\SBRA{\frac{u_1^{1/q_2}}{u_2^{1/q_1}}\PAR{\frac{q_1-q_2}{u_1q_1-u_2q_2}}^{1/q_2-1/q_1}}^{n/2}
\end{equation}
for all functions $f$.
  \end{enumerate}
  \end{ethm}
\begin{eproof}
The proof partly follows the proof of Theorem~\ref{theo-cdn}.  We first assume $(i)$ and prove $(ii)$. For that purpose we consider the map $H$ defined on $\dR^+$ by
$$
H(t)=\PT{u}\PAR{e^{q\Qt f}}^{1/q}
$$
where $u$ and  $q$ are functions of $t \geq 0$. By differentiation and using the Hamilton-Jacobi equation we get  
$$
\frac{q^2}{q'} H' H^{q-1} =\ent{\PT{u}}{e^{q\Qt f}}+u'\frac{q^2}{q'}\PT{u}\PAR{e^{q \Qt f}L \Qt f}+q^2\PAR{\frac{u'q}{q'}-\frac{1}{2q'}}\PT{u}\PAR{\Gamma(\Qt f)e^{q\Qt f}}.
$$ 
Then, by Theorem~\ref{theo-cdn}, assumption $(i)$ implies the logarithmic Sobolev inequality~\eqref{eq-logsobn} which, in the form~\eqref{eq-logsobn'}, writes
$$
\ent{P_u}{e^{q\Qt f}}\leq uq(1-\lambda)\PT{u}\PAR{e^{q \Qt f}L \Qt f}+uq^2\PT{u}\PAR{\Gamma(\Qt f)e^{q\Qt f}}+\frac{n}{2}A_\lambda\PT{u}\PAR{e^{q\Qt f}}
$$
for all $\lambda >0,$ where again $A_{\lambda} = \lambda - 1 - \log \lambda.$ Then we let $q$, $u$ and $\lambda$ solve   
$$
\left\{
\begin{array}{l}
\disp u'q+uq'-uq'\lambda=0\\
\disp q'u+u'q=\frac{1}{2},
\end{array}
\right.
$$
so that
$$
\frac{q^2(t)}{q'(t)} H' H^{q-1} (t) \leq \frac{n}{2}A_{\lambda(t)} \, H^q (t),
$$
and then
$$
\PT{u(t) }\PAR{e^{q(t)\Qt f}}^{1/q(t)}\leq \PT{u(0)}\PAR{e^{q(0)\QT{0}{f}}}^{1/q(0)} \exp \Big( \frac{n}{2} \int_0^t A_{\lambda(r)} \frac{q'(r)}{q(r)^2} \, dr \Big)
$$
if $q$ is increasing. The two equations give $(uq)'=1/2$, that is $u(r)=\PAR{1/2 r+C}/{q}$ for a constant $C$ and 
$$
\lambda (r)=\frac{q(r)}{q'(r)(r+2C)}.
$$
As in the proof of Theorem~\ref{theo-cdn} we minimize the quantity $  \int_0^t A_{\lambda(r)} \frac{q'(r)}{q(r)^2} \, dr$ by letting
$$
q(r)=\frac{r+2C}{\alpha r+\beta}
$$
on $[0,t]$, where $\beta$ is a constant and  
$$
\alpha=\frac{u_1-u_2}{u_1q_1-u_2q_2}
$$
if $q_1 = q(t), q_2 = q(0), u_1 = u(t)$ and $u_2 = u(0).$ We obtain 
$$
\lambda (r) =\frac{1}{1-\alpha q(r)}
$$
so that
$$
\int _0^{t}\frac{q'(r)}{q^2(r)}A_{\lambda(r)} dr=\int _{q_2}^{q_1}\PAR{\frac{1}{1-\alpha q}-1+\log\PAR{1-\alpha q}}\frac{dq}{q^2}
=
\frac{\log\PAR{1-\alpha q_2}}{q_2} - \frac{\log\PAR{1-\alpha q_1}}{q_1}.
$$
 Then the relations $q(r) (\alpha r+\beta) = r+2C = 2 q(r) u(r)$ lead to the condition $t = 2 (u_1 q_1 - u_2 q_2)$, and we obtain $(ii).$
 
 Then we assume $(ii).$ As $ \varepsilon$ goes to 0  a first order Taylor expansion of~\eqref{eq-hj2} with $q_2=1$, $u_2=u_1(1+a  \varepsilon)$ and $t=2 \varepsilon u_1(1-a)$, gives the local logarithmic Sobolev inequality~\eqref{eq-logsobn}, which implies $(i)$. 
\end{eproof}

Let us now see how these hypercontractive bounds on Hamilton-Jacobi equations translate into transportation inequalities.  
 %%%%%%%%%%%%%%%%%%%%%%%%%%%%
 \subsection{Application to transportation inequalities}
 \label{sec-transine}
 %%%%%%%%%%%%%%%%%%%%%%%%%%%%%
Let $(M,d)$ be a complete Riemannian manifold. Then the Wasserstein distance between two probability measures $\mu$ and $\nu$ on $M$ is defined by
$$
W_2(\nu,\mu)=\inf_{\pi} \sqrt{\int\frac{d(x,y)^2}{2}d\pi(x,y)}
$$
where the infimum is taken over all probability measures $\pi$ on $M\times M$ such that$$
\int( f(x)+g(y))d\pi(x,y)=\int f \, d\nu+\int g \, d\mu 
$$ 
for all bounded and measurable  functions $f$ and $g$ on $M.$

\medskip

The Wasserstein distance is linked with logarithmic Sobolev inequalities and hypercontractivity bounds on the Hamilton-Jacobi semigroups by the following argument of~\cite{bgl}. Let $C >0$ and a probability measure $\mu$ on $M$ satisfy the logarithmic Sobolev inequality
\begin{equation*}
%\label{eq-islf}
\ent{\mu}{f^2}\leq {2C}\int |\nabla f|^2d\mu
\end{equation*}
for all maps $f$ on $M,$ where $\vert \nabla f \vert$ stands for the Riemannian length of the gradient of $f$.   Then, for $t =1$ and $a$ going to $0$, the hypercontractive bound~\eqref{eq-bgl} leads to
\begin{equation}
\label{eq-bg}
\int \exp\PAR{\frac{1}{C}{\QT{1}{f}}}d\mu\leq\exp\PAR{\frac{1}{C}\int fd\mu}
\end{equation}
for all Lipschitz functions $f.$ Now, following S. Bobkov and F. G\"otze~\cite{bobkov-gotze}, this is equivalent to the inequality
$$
\int \QT{1}f h d\mu-\int h d\nu  \leq C \, \ent{\mu}{f}
$$
for all Lipschitz functions $f$ and all probability densities $h$ with respect to $\mu$, by the variational formulation of the entropy. But, by the
Kantorovich duality formulation,
$$
W_2^2(h \mu,\mu)=\sup_{f, g} \BRA{\int g h d\mu-\int f d\mu}, 
$$
where the supremum is taken over all $f$ and $g$ such that
$$
g(x) + f(y) \leq \frac{1}{2} d(x,y)^2
$$
for all $x$ and $y$ (see~\cite[Chapters 5 and 22]{villani09}). In fact
\begin{equation}\label{eq-Kanto}
W_2^2(h \mu,\mu) =\ \sup_{f} \BRA{\int \QT{1}f h d\mu-\int f d\mu} 
\end{equation}
since, by~\eqref{eq-hjma},  
$$
\QT{1}f(x)=\inf_{x\in M}\BRA{f(y)+\frac{1}{2}d(x,y)^2}
$$
so is optimal among all $g$, for fixed $f.$ In the end, the measure $\mu$ satisfies the transportation or Talagrand  inequality
\begin{equation}
\label{eq-tala-def}
W_2^2(h\mu,\mu)\leq C \, \ent{\mu}{h}
\end{equation}
for all probability densities $h$ with respect to $\mu.$ This inequality was introduced by M. Talagrand in~\cite{tala} for the Gaussian measure. Its implication by a logarithmic Sobolev inequality is called the Otto-Villani Theorem~\cite{otto-villani} and has first been proved differently  (see also~\cite{gozlan-leonard} and~\cite{villani09}). 

\medskip

We now again turn to the local level, and investigate how the curvature-dimension criterion is equivalent to local transportation inequalities, which correspond to the local logarithmic Sobolev inequalities~\eqref{eq-localsob} and~\eqref{eq-logsobn} and hypercontractive bounds~\eqref{eq-hjrho} and~\eqref{eq-hj2}.

%Since  the hypercontractive bound obtained in Theorem~\ref{theo-hjcdn} is stronger than~\eqref{eq-bgl}, the corresponding Talagrand inequality obtained by using~\eqref{eq-hj2} will be stronger than the classical inequality.

%%%%%%%%%%%%%%
\begin{ethm}\label{theo-hjrhot}
Let $(\PT{t})_{t \geq 0}$ be a diffusion Markov semigroup on $M$ such that $\Gamma(f)=\ABS{\nabla f}^2$ for all functions $f$, where $\vert \nabla f \vert$ stands for the Riemannian length of the gradient of $f$.  

\begin{itemize}
\item Given $\rho\in\dR$  the following assertions are equivalent :
\begin{enumerate}[(i)]
\item
 the semigroup $(\PT{t})_{t \geq 0}$ satisfies the ${CD}(\rho,\infty)$ criterion;
 \item 
for all $u \geq 0$, all $x\in M$ and all $h\geq0$ such that $\PT{u} h (x)=1$,
\begin{equation}
\label{eq-tala-rho}
W_2^2(h\PT{u}^x,\PT{u}^x)\leq \frac{1-e^{-2\rho u}}{\rho}\ent{\PT{u}^{x}} {h},
\end{equation}
where $h\PT{u}^x$ is the probability measure defined by  $\int \phi \,d(h\PT{u}^x) =\PT{u}(\phi h)(x)$ for all functions $\phi$.
% {\red l'ecriture est assez nulle, il faut peut-etre l'ecrire avec les operateurs $\Pt^*$}
  \end{enumerate}
  If $\rho=0$ then $\PAR{1-e^{-2\rho u}}/{\rho}$ will be replaced by  $2u$.
\item Given $n\geq 1$ the following assertions are equivalent :
\begin{enumerate}[(i)]
\item
 the semigroup $(\PT{t})_{t \geq 0}$ satisfies the ${CD}(0,n)$ criterion;
\item
for all $u_2,u_1\geq 0$, all $x\in M$ and all $h\geq0$ such that $\PT{u_1} h (x)=1$, 
\begin{equation}
\label{eq-talann}
W_2^2(h\PT{u_1}^x,\PT{u_2}^x)\leq 2 u_1 \PAR{\ent{{\PT{u_1}^x}}{h}+\frac{n}{2} A_{u_2/u_1}}.
\end{equation}
  \end{enumerate}
%  {\red fb: il semble que le $2 u_1$ soit en facteur dans le membre de droite, a verifier IG: OK}
  \end{itemize}
  \end{ethm}
%%%%%%%%%%%%

\begin{eproof}
We first prove  the first part of the theorem.  Let us first assume that the ${CD}(\rho,\infty)$ criterion holds, so that inequality~\eqref{eq-hjrho} is satisfied by Theorem~\ref{theo-hjrho}.  We let $t=1$ and $u > 0$ fixed, and take the limit in~\eqref{eq-hjrho} as $q_2$ goes to $0$. We obtain 
\begin{equation}
 \label{eq-bg2}
 \PT{u}\PAR{e^{q_1\QT{1} f}}\leq \exp\PAR{q_1\PT{u}{f}}
 \end{equation}
with $q_1={\rho}/\PAR{1-e^{-2\rho u}}$.
This is equivalent to~\eqref{eq-tala-rho} by the dual representation~\eqref{eq-bg} of S. Bobkov and F. G\"otze, applied to the measure $\PT{u}^x$, which proves $(ii).$

Then we assume $(ii)$, that is, \eqref{eq-bg2} for all $f$. We let $f = 1+ \varepsilon g$ and perform a second order Taylor expansion for $ \varepsilon$ going to 0. We obtain the local Poincar\'e inequality
$$
\var{\PT{u}}{g}\leq \frac{1-e^{-2\rho u}}{\rho}\PT{u}\PAR{\Gamma(g)},
$$
which is equivalent to the $CD(\rho,\infty)$ criterion (see~\cite[Proposition 3.3]{bakrytata} for instance).

\medskip
Then we prove  the second part of the theorem.  Let us first assume that the ${CD}(0,n)$ criterion holds, so that inequality~\eqref{eq-hj2} is satisfied by Theorem~\ref{theo-hjcdn}. We let $x \in M$ and $u_1 >0$ be fixed and we apply the inequality~\eqref{eq-hj2} with $s=0$, $t=1$, $q_1 = (2 u_2 q_2 +1)/(2 u_1)$ and $q_2$ going to $0$. We obtain 
\begin{equation}
\label{eq-bg3}
\PT{u_1}\PAR{e^{q_1\QT{1}{f}}}\leq e^{q_1\PT{u_2}f}\exp\PAR{\frac{n}{2}A_{u_2/u_1}},
\end{equation}
where $q_1=1/(2u_1)$, that is,
$$
\int \exp \Big( q_1\QT{1}{f}(y) - q_1 \PT{u_2}f (x) - \frac{n}{2}A_{u_2/u_1} \Big) \PT{u_1}^x (dy) \leq 1.
$$
Hence, by the variational formulation of the entropy,
$$
q_1 \Big( \int \QT{1}{f} (y) \, h (y) \, \PT{u_1}^x (dy) - \int f(y) \, \PT{u_2}^x (dy) \Big) \leq \ent{\PT{u_1}}{h}(x)+\frac{n}{2} A_{u_2/u_1}
$$
for all fixed $x$ and all positive functions $h$ such that $\PT{u_1}{h} (x) = \int h d\PT{u_1}^x = 1.$ In other words
$$
q_1 \, W_2^2(h\PT{u_1}^x,\PT{u_2}^x)\leq \ent{\PT{u_1}}{h}(x)+\frac{n}{2} A_{u_2/u_1}
$$
by the Kantorovich duality formulation~\eqref{eq-Kanto}, which proves the inequality~\eqref{eq-talann}.

Then we assume $(ii).$  We first note that $(ii)$ is equivalent to having~\eqref{eq-bg3} for all $f$. Then, as in the proof of Theorem~\ref{theo-cdn}, we perform a second order Taylor expansion of inequality~\eqref{eq-bg3}. We let $f= \varepsilon g$, $u_1=t$, $u_2=t(1+a \varepsilon)$ where $a\in\dR$.  We note in particular that
$$
\QT{1}{ \varepsilon g}= \varepsilon g-\frac{ \varepsilon^2}{2}\Gamma(g)+o( \varepsilon^2)
$$
as $ \varepsilon$ goes to 0. We obtain the bound
$$
\var{\Pt}{g}\leq 2t\Pt(\Gamma(g)) + 4 a t^2 L\Pt g + 2 n t^2 a^2 = 
2t\Pt(\Gamma(g))-\frac{2t^2}{n}\PAR{L\Pt g}^2
$$
for $a = - L \Pt g / n.$ This finally implies the $CD(0,n)$ criterion by~\cite[Proof of Theorem 1]{bakry-ledoux06}.
\end{eproof}
 
\medskip

In Remark~\ref{rem-gross} we saw how the local hypercontractive bound of Theorem~\ref{theo-rho} leads to the Gross Theorem under a $CD(\rho, \infty)$ condition with $\rho >0.$ In the same way, if the semigroup $(\Pt)_{t\geq0}$ satisfies the $CD(\rho,\infty)$ condition with $\rho >0$
and is ergodic in $L^p(d\mu)$ for the probability measure $\mu$, then, as $u$ goes to infinity,  the inequality~\eqref{eq-tala-rho} in the first part of Theorem~\ref{theo-hjrhot} implies the Talagrand inequality~\eqref{eq-tala-def} for $\mu$, with $C = 1/\rho.$
 
 Then we saw in Remark~\ref{rem-ou}  how to refine the Nelson Theorem on the Ornstein-Uhlenbeck semigroup by using a $CD(0,n)$ condition on the heat semigroup. In the same way we now refine the usual Talagrand inequality for the Gaussian measure
as a simple consequence of the second part of Theorem~\ref{theo-hjrhot}, in $\dR^n$ equipped with                     the Euclidean norm:

\begin{ecor}
The standard Gaussian measure on $\dR^n$ satisfies the refined Talagrand inequality 
\begin{equation}
\label{eq-talan}
W_2^2(h\gamma,\gamma)\leq \frac{1}{2}\int \Delta h \, d\gamma+{n}\PAR{1-\exp\PAR{{\frac{1}{2n}\int \Delta h \, d\gamma-\frac{1}{n}\ent{\gamma}{h}}}}
\end{equation}
for all probability densities $h$ with respect to $\gamma$. This inequality implies the classical Talagrand inequality 
$$
W_2^2(h\gamma,\gamma)\leq \ent{\gamma}{h}
$$
and the refined Poincar\'e inequality 
$$
\var{\gamma}{f}\leq \int |\nabla f|^2d\gamma-\frac{1}{2n}\PAR{\int \Delta fd\gamma}^2
$$ 
for all smooth functions $f$.
\end{ecor}
\begin{eproof}
Let indeed $L$ be the Laplacian operator and $(\Pt)_{t\geq 0}$ be the heat semigroup, which again satisfies the $CD(0,n)$ condition.  We apply $(ii)$ in the second part of  Theorem~\ref{theo-hjrhot} with $x=0$ and $u_2 = 1/2$, in the form
$$
W_2^2(h \PT{u}^0, \gamma) \leq 2 u \big( \ent{\PT{u}^0}{h}+\frac{n}{2} A_{1/(2u)} \big)
$$
for all $u >0.$ 

Then we consider a map $g$ on ${\mathbb R}^n$ such that $\displaystyle \int  g \, d\gamma =1$. We let 
$$
h (x)= g(x) \frac{ \gamma(x)}{\PT{u}^0 (x)},
$$
so that $h \, \PT{u}^0  = g \, \gamma$ has mass $1$ and
$$
W_2^2(g \gamma, \gamma) \leq  2 u \big( \ent{\PT{u}^0}{h}+\frac{n}{2} A_{1/(2u)} \big).
$$
We note that
$$
\ent{\PT{u}^0}{h} = \ent{\gamma}{g} + \frac{n}{2} \log (2u) + \big(  \frac{1}{4u} -  \frac{1}{2} \big) \int \vert x \vert^2 g(x) \, d\gamma(x)
$$
where
$$
\int \vert x \vert^2 g(x) \, d\gamma(x) = - \int x  g(x) \cdot \nabla \gamma(x) \, dx = \int \nabla \cdot (x g(x)) \gamma(x) \, dx = n +  \int \Delta g \,d \gamma .
$$
We finally obtain the bound
\begin{equation}\label{eq-talan-mu}
W_2^2(g\gamma,\gamma)\leq \lambda \, \ent{\gamma}{g} + n(1-\lambda+ \lambda\log \lambda)+\frac{1}{2}(1-\lambda)\int \Delta g \, d\gamma
\end{equation}
for all positive $\lambda \, ( = 2 u).$ 

The choice $\lambda =1$ leads to the usual Talagrand inequality
\begin{equation}\label{eq-tala-gauss}
W_2^2(g \gamma,\gamma)\leq \ent{\gamma}{g}
\end{equation}
for the Gaussian measure, whereas the optimal choice leads to the bound~\eqref{eq-talan}.  

We may recover that the inequality~\eqref{eq-talan} implies the usual bound~\eqref{eq-tala-gauss} by applying the inequality $e^x\geq x+1$ for $x\in\dR.$

\medskip

We now prove that the refined Talagrand inequality~\eqref{eq-talan} implies the refined Poincar\'e inequality
$$
\var{\gamma}{f}\leq \int|\nabla f|^2d\gamma-\frac{1}{2n}\PAR{\int \Delta fd\gamma}^2.
$$
For that purpose we first note that~\eqref{eq-talan} is in fact equivalent to the bounds~\eqref{eq-talan-mu} for all $\lambda >0,$ so by the Kantorovich duality formulation~\eqref{eq-Kanto} implies the bound
$$
\int \QT{1} f g d\gamma - \int f \, d\gamma \leq \lambda \, \ent{\gamma}{g} + n(1-\lambda+ \lambda\log \lambda)+\frac{1}{2}(1-\lambda)\int \Delta g \, d\gamma
$$
for all maps $f, g$ and $\lambda >0.$ 

Then we consider a map $h$ on $\mathbb R^n$ and let $f(x) =  \varepsilon h(x), g(x) = 1 +  \varepsilon (h(x) + a (\vert x \vert^2 - n))$ and $\lambda = 1 + 2 a  \varepsilon.$ We perform a second order Taylor expansion for $ \varepsilon$ going to $0$ and obtain the bound
$$
\var{\gamma}{h}\leq \int|\nabla f|^2d\gamma - a \int \Delta h \, d\gamma + \frac{n a^2}{2}
$$ 
which leads to the refined Poincar\'e inequality after optimization over $a.$
\end{eproof}

\begin{erem}
For non smooth probability densities $h$ with respect to $\gamma$, the inequality~\eqref{eq-talan} writes
$$
W_2^2(h\gamma,\gamma)\leq \frac{1}{2} \int \vert x \vert^2 h(x) \, d\gamma(x) - \frac{n}{2}  +{n}\PAR{1-\exp\PAR{{\frac{1}{2n}\int \vert x \vert^2 h(x) \, d\gamma(x) - \frac{1}{2} - \frac{1}{n}\ent{\gamma}{h}}}},
$$
which can be useful in concentration of measure arguments.

\end{erem}
%\begin{erem}
%The inequality~\eqref{eq-talan} can be proved, using the refined logarithmic Sobolev inequality proved in \cite[Proposition 2]{bakry-ledoux06}, for all nonnegative function $f$, such that $\int fd\gamma=1$
%$$
%\ent{\gamma}{f}\leq \frac{1}{2}\int \Delta fd\gamma+\frac{n}{2}\log\PAR{1-\frac{1}{n}\int f\Delta\PAR{\log f}d\gamma}.
%$$
%Then  the refined Talagrand inequality is between the refined Poincar\'e and the refined logarithmic Sobolev inequalities. 
%{\red a verifier}
%\end{erem}

%%%%%%%%%%%%%%%%%%%%%%%%%%%%%%%%%%%%%%%%%%%%%%%%%%%%%%

\section{Hypercontractivity for L\'evy operators}
\label{sec-levy}
%%%%%%%%%%%%%%%%%%%%%%%%%%%%%%%%%%%%%%%%%%%%%%%%%%%%%%
\subsection{The general case}
A L\'evy operator is by definition the infinitesimal generator of a L\'evy process. 

For a  pure jump process, the L\'evy-Khinchine formula asserts that there exists a positive measure $\nu$ on $\dR$ such that $\nu(\BRA{0})=0$ and  
$$
\int_{\dR}\min\PAR{1,|z|^{2}}d\nu(z)<\infty,
$$
and such  that its generator  is given by  
\begin{equation}
\label{eq-defle}
\mathcal I (f) (x) = \int \big( f(x+z) -f (x) -\nabla f (x) \cdot z h(z) \big) d\nu(z)
\end{equation}
for all functions $f$, where $h(z)=1/(1+|z|^2)$. We let $(\Lt)_{t\geq 0}$ be the Markov semigroup associated to  $\mathcal I$. Classical results on L\'evy semigroups are given in~\cite{applebaum}.

In our context, a L\'evy semigroup is a non-diffusive Markov semigroup satisfying  a $CD(0,\infty)$ criterion, see \cite{bakryemery}.  It also satisfies the local logarithmic Sobolev inequality
\begin{equation}
\label{eq-logsoblevy}
\ent{\Lt}{ f} \leq t \, \Lt \PAR{\int D (f(\cdot +z), f(\cdot)) d\nu(z)}
 \end{equation}
 for all positive functions $f$,  where, for all positive $u$ and $v$,
 $$
D (u,v) = u \log \frac{u}{v} - (u-v).
$$
The function $D$ is known as the Bregman distance associated to the function $x\mapsto x\log x$. Inequality~\eqref{eq-logsoblevy} has been proved by C. An\'e and M. Ledoux in~\cite{ane-ledoux} and by L. Wu in~\cite{wu00}, one can also see~\cite{chafai04} for the more general setting of $\Phi$-entropy inequalities. 

A L\'evy operator is not a diffusion operator, so we cannot directly apply Theorem \ref{theo-rho} with $\rho =0.$  Nevertheless we obtain the following weaker hypercontractive bound :
 \begin{ethm}
\label{thm-levy}
For all $0< s \leq t$ and $0<q_2 \leq q_1<1$ such that
\begin{equation}
\label{eq-levycond}
{q_2-1} =\frac{q_2}{q_1}\frac{t}{s} \PAR{q_1-1},
\end{equation}
then
 \begin{equation}
 \label{eq-levyrho}
 \LT{t}{(f^{q_1})}^{1/q_1} \leq \LT{s}{(\LT{t-s}{(f)}}^{q_2})^{1/q_2}
 \end{equation}
for all positive functions $f$. If now $q_2-1=\frac{t}{s}(q_1-1)$ and if moreover $\nu\PAR{\dR}=C<\infty$, then
\begin{equation}
\label{eq-levyc}
\LT{t}{(f^{q_1})}^{1/{q_1}} \leq \LT{s}{(\LT{t-s}{(f)}}^{q_2})^{1/q_2} \exp \Big[ C \frac{s(t-s)(1-q_2)^2}{q_2(sq_2+t-s)} \Big].
\end{equation}
 \end{ethm}
 \begin{erem}
In the case when $\nu(\dR)<\infty$ then the two hypercontractive bounds in Theorem~\ref{thm-levy} correspond to two extremal points of view: inequality~\eqref{eq-levyrho} gives a finer bound whereas inequality~\eqref{eq-levyc} gives a larger range of parameters $q_2$ for given $q_1.$  But one can obtain  interpolating bounds for intermediate parameters $q_2$.  

Given $s, t$ and $q_1$, the bound~\eqref{eq-levyrho} is obtained by a parameter $q_2$ closer to $q_1$ than  in the bound~\eqref{eq-majrho2} obtained in the diffusive case under a $CD(0, \infty)$ condition. To reach the same gain of integrability we need to impose a finite mass assumption on $\nu$ and an extra constant in the bound~\eqref{eq-levyc}. Also the parameters $q_1$ and $q_2$ are restricted to $(0,1).$ This is the price to pay when passing from the diffusive setting of Theorem~\ref{theo-rho} to this non diffusive setting.
\end{erem}

\begin{eproof}
Let, for $s\in[0,t]$,
$$
\psi(s)=\LT{s}{\PAR{\LT{t-s}{ \PAR{f}^{q}}}}^{1/q},
$$  
where $q$ depends on $s$. By differentiation
 $$
 \psi'(s) \frac{q^2}{q'} \psi(s)^{q-1}=  { \ent{\LT{s}}{ g^q} + \frac{q}{q'} \LT{s}{\PAR{{\mathcal I}(g^q) - q g^{q-1} I(g)}}}
 $$
for all $s\in[0,t],$ where $g={\LT{t-s}{ \PAR{f}}}$.  The logarithmic Sobolev inequality~\eqref{eq-logsoblevy} applied to the function $g^q$ implies
$$
\psi'(s) \frac{q^2}{q'} \psi(s)^{q-1} \leq  \LT{s}{\PAR{\int \Phi(g(\cdot+z),g(\cdot)) \, d\nu(z)}} 
 $$ 
 where 
$$
\Phi(Z,X)={{sZ^q\log\frac{Z^q}{X^q}-s(Z^q-X^q)+\frac{q}{q'}(Z^q-X^q)-\frac{q^2}{q'}X^{q-1}(Z-X)}}.
$$
Let now $q$ be an increasing function on $[0,t]$ such that $q\in(0,1)$ and  
$$
 s + q \frac{q-1}{q'} =0,
 $$
for instance such that
$$
\frac{q(t) - 1}{q(s) - 1} = \frac{q(t)}{q(s)} \frac{s}{t}.
$$ 
Then $\Phi(Z,X)\leq 0$ for all $X,Z\geq 0$, so that $\psi'(s) \leq 0$ and finally $\psi(t) \leq \psi(s).$ This is 
inequality~\eqref{eq-levyrho} provided $q_1 = q(t)$ and $q_2 = p(t)$ satisfy~\eqref{eq-levycond}.

 If now $q$ is increasing and satisfies
$$
 s +  \frac{q-1}{q'} =0,
 $$
then
$\Phi(Z,X)\leq \frac{(q-1)^2}{q'}X^q$, which implies
$$
\psi'(s) \frac{q(s)^2}{q'(s)} \psi(s)^{q-1}\leq C\psi(s)^{q}
 $$ 
where $C = \int d\nu$ is assumed to be finite. Integrating this inequality implies~\eqref{eq-levyc}. 
\end{eproof}

\subsection{The case of $\alpha$-stable L\'evy processes}
The semigroup given by a $\alpha$-stable L\'evy process can be represented  as 
$$
\LT{t}^\alpha f(x)=K_t^\alpha\star f(x)=\int K_t^\alpha(x-y) f(y)dy
$$
for all $t \geq 0$, all $x\in\dR^n$ and all smooth functions $f$ on $\dR^n$, where 
$$
K_t^\alpha(x)=\frac{1}{(2\pi)^{n/2}}\mathcal F\PAR{e^{-t\ABS{\cdot}^\alpha}}(x);
$$
here $\mathcal F (g)$ is the Fourier transform of  a function $g\in L^1(dx)$, given by 
$$
\mathcal F(g)(x)=\frac{1}{(2\pi)^{n/2}}\int e^{-iu\cdot x}f(u)du.
$$
The generator of such a semigroup will be denoted $\mathcal I^\alpha$. It is given by~\eqref{eq-defle} with the L\'evy measure
$$
d\nu_\alpha(x)=\frac{c_\alpha}{\ABS{x}^{n+\alpha}}dx
$$
for a positive constant $c_{\alpha}$. It is  the fractional Laplacian operator $-(-\Delta)^{2/\alpha}$, which for $\alpha =2$ is the Laplacian. All details can be found in~\cite{applebaum}. 

 This generator is a non-diffusive operator but has the Lebesgue measure as a reversible measure and using a fractional Euclidean logarithmic Sobolev inequality one can obtain the following hypercontractive bound :  there exists a constant $\mathcal A$ such that  
\begin{equation} \label{eq-diff1}
\NRM{\LT{t}^\alpha f}_{L^{q_2}(dx)} \leq \NRM{ f}_{L^{q_1}(dx)} \left( \frac{{\mathcal A}n (q_2 - q_1)}{2
  \alpha t}\right)^{\frac{n (q_2 -q_1)}{\alpha q_1q_2}}
  \frac{{q_1}^{n/(q_2\alpha)}}{{q_2}^{n/(q_1\alpha )}}
\end{equation}
for all $q_2>q_1 \geq 2$ and all positive functions $f$ (see~\cite{gentil-imbert09} for the proof).

\medskip

Let us now see, following Remark~\ref{rem-ou}, how Theorem~\ref{thm-levy} leads to a hypercontractive bound for the L\'evy-Ornstein-Uhlenbeck semigroup, which is an Ornstein-Uhlenbeck process driven by a L\'evy process instead of a Brownian motion. 

\medskip

As for the heat equation, 
$$
\LT{t}^\alpha(f)(0)=\NRM{\TT{\frac{2}{\alpha}\log\PAR{\alpha t}}f}_{L^1(K^{\alpha}_{1/\alpha})},
$$
so that if  $0<q_2 \leq q_1<1$ satisfy~\eqref{eq-levycond} then  inequality~\eqref{eq-levyrho} at $x=0$ gives  
$$
\NRM{\TT{\frac{2}{\alpha}\log (\alpha t)} {f}}_{L^{q_1}(K_{1/\alpha})}\leq \NRM{\TT{\frac{2}{\alpha}\log (\alpha s)} {\LT{t-s}^\alpha f}}_{L^{q_2}(K^{\alpha}_{1/\alpha})},
$$
%{\red verifier les p et q}
that is, 
$$
\NRM{{g}}_{L^{q_1}(K_{1/\alpha}^\alpha)}\leq \NRM{\TT{\frac{2}{{\alpha}}\log (\alpha s)} {\LT{t-s}^\alpha\TT{-\frac{2}{\alpha}\log (\alpha t)} {g }}}_{L^{q_2}(K^{\alpha}_{1/\alpha})}
$$
for all positive functions $g$. But $\LT{t}^\alpha$ and $\Tt$ commute according to
$$
\LT{t}^\alpha \TT{a}=\TT{a}\LT{te^{a\alpha/2}}^\alpha,
$$
so that
$$
\NRM{{g}}_{L^{q_1}(K_{1/\alpha}^\alpha)}\leq 
\NRM{\TT{\frac{2}{\alpha}\log \frac{s}{t}} {\LT{\frac{1}{\alpha}\PAR{1-\frac{s}{t}}}^\alpha {g}}}_{L^{q_2}(K^{\alpha}_{1/\alpha})}
 = 
\NRM{\TT{- 2 t } {\LT{\frac{1}{\alpha}\PAR{1- e^{ - \alpha t}}}^\alpha {g}}}_{L^{q_2}(K^{\alpha}_{1/\alpha})}
 $$
for $s/t=e^{-\alpha t}.$ But the semigroup 
$$
\PT{t}^{LOU}=\TT{-2t} {\LT{\frac{1}{\alpha}\PAR{1-e^{- \alpha t}}}^\alpha },
$$
is the L\'evy-Ornstein-Uhlenbeck semigroup, solution of
$$
{\partial_ t}\PT{t}^{LOU}(g) = \mathcal I^\alpha(\PT{t}^{LOU}(g))-x\cdot\nabla\PT{t}^{LOU}(g).
$$
Hence we have obtained the following result :
\begin{ecor}
The Markov semigroup of generator  $L=\mathcal I^\alpha-x\cdot\nabla$ satisfies the hypercontractive bound 
$$
\NRM{{f}}_{L^{q_1}(K_{1/\alpha}^\alpha)}\leq \NRM{ \PT{t}^{LOU}{f}}_{L^{q_2}(K^{\alpha}_{1/\alpha})},
$$ 
for all $t \geq 0$ and $0 < q_2\leq q_1< 1$ such that 
$$
{q_2-1}=\frac{q_2}{q_1}\PAR{q_1-1}e^{\alpha t}.
$$
\end{ecor}
\begin{erem}
This is a kind of reverse hypercontractive bound on the L\'evy-Ornstein-Uhlenbeck semigroup. This bound is classical for the usual Ornstein-Uhlenbeck semigroup, see~\eqref{eq-derder}.

In our case, parameters $q_1,q_2$ are closer than in the classical case for two reasons : our semigroup is not a diffusion in the sense of the definition~\eqref{eq-difff} and an $\alpha$-stable process is less `{\it diffusive}" than the Brownian motion. 
 
Note that the probability measure $K^{\alpha}_{1/\alpha}$ is invariant for the generator $\mathcal I^\alpha-x\cdot\nabla$ but it is not reversible. To our knowledge this is the first instance of a hypercontractive semigroup which is not a diffusion and has a non reversible invariant measure. 
\end{erem}

{\bf Acknowledgements.} This research was supported in part by the ANR project EVOL. The third author thanks the members of UMPA at the Ecole Normale Sup\'erieure de Lyon for their kind hospitality.

\noindent

\medskip\noindent

\noindent
Institut de Math\'ematiques de Toulouse\\
Universit\'e de Toulouse\\
31062 Toulouse - France\\
bakry@math.univ-toulouse.fr 

\medskip

\noindent
Ceremade\\
Universit\'e Paris-Dauphine\\
Place du Mar\'echal De Lattre De Tassigny\\
75116 Paris - France\\
bolley@ceremade.dauphine.fr\\
gentil@ceremade.dauphine.fr\\

\end{document}